\documentclass{amsart}
\usepackage{amssymb}
\usepackage{mathrsfs}
\usepackage{cases}
\usepackage{amsmath}
\usepackage{amsfonts}
\usepackage{pifont}
\usepackage{graphicx,amssymb,mathrsfs,amsmath}
\usepackage{tocvsec2}
\newtheorem{thm}{Theorem}[section]

\newtheorem{prop}[thm]{Proposition}
\theoremstyle{definition}
\newtheorem{defn}[thm]{Definition}
\newtheorem{example}[thm]{Example}
\theoremstyle{remark}
\newtheorem{rem}[thm]{Remark}
\numberwithin{equation}{section}

\begin{document}

\title{${\mathcal F}$-Hypercyclic and disjoint ${\mathcal F}$-hypercyclic properties of binary relations over topological spaces}

\author{Marko Kosti\' c}
\address{Faculty of Technical Sciences,
University of Novi Sad,
Trg D. Obradovi\' ca 6, 21125 Novi Sad, Serbia}
\email{marco.s@verat.net}

\begin{abstract}
In this  paper, we examine various types of ${\mathcal F}$-hypercyclic (${\mathcal F}$-topologically transitive) and disjoint ${\mathcal F}$-hypercyclic (disjoint ${\mathcal F}$-topologically transitive) properties 
of binary relations over topological spaces. We pay special attention to finite structures like simple graphs, digraphs and tournaments, providing a great number of illustrative examples.
\end{abstract}

{\renewcommand{\thefootnote}{} \footnote{2010 {\it Mathematics
Subject Classification. 47A16, 47B37, 47D06.
\\ \text{  }  \ \    {\it Key words and phrases.}  ${\mathcal F}$-hypercyclic binary relations, ${\mathcal F}$-topologically transitive binary relations, disjoint ${\mathcal F}$-hypercyclic binary relations, disjoint ${\mathcal F}$-topologically transitive binary relations, digraphs.
\\  \text{  } The author is partially supported by grant 174024 of Ministry
of Science and Technological Development, Republic of Serbia.}}

\maketitle

\section{Introduction and preliminaries}\label{into}

The notion of a continuous linear frequently hypercyclic operator acting on a separable Fr\' echet space was introduced by F. Bayart and S. Grivaux in \cite{bay1} (2006). From then on, a great number of authors working in the field of linear topological dynamics has analyzed the notion of frequent hypercyclicity, various generalizations of this concept and certain applications to abstract differential equations. Recently, upper frequent hypercyclic linear operators and ${\mathcal F}$-transitive linear operators have been investigated by A. Bonilla, K.-G. Grosse-Erdmann  \cite{boni-upper} and J. B\`es, Q. Menet, A. Peris, Y. Puig \cite{biba}.  For more details on the subject, we refer the reader to \cite{Bay}-\cite{biba}, \cite{boni}-\cite{boni-ERAT}, \cite{Grosse}, \cite{menet} and references cited therein.

The main aim of this paper is to continue the research studies \cite{marina}-\cite{marina-prim} and \cite{first-part}-\cite{extensions}. We analyze ${\mathcal F}$-hypercyclic (${\mathcal F}$-topologically transitive) and disjoint ${\mathcal F}$-hypercyclic (disjoint ${\mathcal F}$-topologically transitive) properties of binary relations over topological spaces, focusing special attention to finite topological spaces which do not have a linear vector structure. Concerning similar problematic, one may refer e.g. to the papers by R. A. Mart\' inez-Avendano \cite{martinez}, where the author has investigated hypercyclic shifts on weighted $L^{p}$ spaces of directed trees, P. Namayanja \cite{south-africa}, where 
chaotic phenomena in a transport equation on a network have been studied with the use of adjacency matrices of infinite line graphs, and C.-C. Chen \cite{chenuljica}, where the author has investigated hypercyclic and chaotic operators on $l^{p}$ spaces of Cayley
graphs. We present plenty of results and illustrative examples for simple graphs, digraphs and tournaments. With the exception of paper \cite{extensions}, where we have recently analyzed ${\mathcal F}$-hypercyclic extensions and disjoint ${\mathcal F}$-hypercyclic extensions of binary relations over topological spaces, the notions of ${\mathcal F}$-hypercyclicity and ${\mathcal F}$-topological transitivity have not been considered elsewhere in such a general framework.
And, more to the point, with the exception of paper \cite{extensions}, disjoint ${\mathcal F}$-hypercyclicity and disjoint ${\mathcal F}$-topological transitivity seem to be not considered elsewhere even for linear continuous operators acting on Banach spaces. 

The organization and main ideas of paper are briefly described as follows. In Section \ref{into}, we introduce the notions of ${\mathcal F}$-hypercyclicity and disjoint ${\mathcal F}$-hypercyclicity for binary relations, giving also a few noteworthy observations and elementary consequences of definitions. We divide the third section of paper in three separate subsections. In Subsection 3.1, we analyze ${\mathcal F}$-hypercyclicity and disjoint ${\mathcal F}$-hypercyclicity for general binary relations, on finite or infinite topological spaces, having or not a certain number of loops. We slightly extend the implications (a) $\Rightarrow$ (b) $\Rightarrow$ (c) $\Rightarrow$ (d) of a recent result by A. Bonilla, K.-G. Grosse-Erdmann  \cite[Theorem 15]{boni-upper} in this context, and reformulate the notions introduced in the second section in terms of appropriate conditions on adjacency matrices (it is worth noting that \cite[Theorem 15]{boni-upper} is exceptional in the existing theory of linear topological dynamics because it is a rare result in which the pivot spaces do not need to be equipped with linear vector structures). In Subsection 3.2, we focus our attention to the simple graphs. In Proposition \ref{poka}, we firstly prove that the notions of $d{\mathcal F}$-hypercyclicity ($d{\mathcal F}$-topological transitivity) and strong $d{\mathcal F}$-hypercyclicity (strong $d{\mathcal F}$-topological transitivity)
coincide in the case that ${\mathcal F}=P(P({\mathbb N})) \setminus \{\emptyset \},$ which is unquestionably the best explored in the existing literature.
For a simple graph $G,$ we introduce the index ${\mathrm S}_{G}$ and give some upper bounds for $S_{G}$ in Theorem \ref{moguce}, concerning connected bipartite graphs, and Theorem \ref{reza}, concerning connected non-bipartite graphs. In Theorem \ref{radio}, we prove that connected non-bipartite graphs $G_{1}, G_{2},\cdot \cdot \cdot, G_{N}$ are always (strongly) $d{\mathcal F}$-hypercyclic (strongly $d{\mathcal F}$-topologically transitive). Excluding Proposition \ref{idiot}, almost all structural results from Subsection 3.3 is devoted to the study of case ${\mathcal F}=P(P({\mathbb N})) \setminus \{\emptyset \}.$ We pay a special attention to the question whether the $d{\mathcal F}$-hypercyclicity ($d{\mathcal F}$-topological transitivity) of a given digraph $G$ (digraphs $G_{1}, G_{2},\cdot \cdot \cdot, G_{N}$) automatically implies the strong $d{\mathcal F}$-hypercyclicity (strong $d{\mathcal F}$-topological transitivity) of $G$ ($G_{1}, G_{2},\cdot \cdot \cdot, G_{N}$). In Proposition \ref{pende-primp}, we prove that, if 
the number of nodes of a digraph $G$ is less than or equal to $4,$ and $G$ is equipped with arbitrary topology, then the ${\mathcal F}$-hypercyclicity of $G$ always implies the strong ${\mathcal F}$-hypercyclicity of $G$. For disjointness, we prove that the $d{\mathcal F}$-hypercyclicity of $G_{1},G_{2},\cdot \cdot \cdot, G_{N}$ always implies the strong $d{\mathcal F}$-hypercyclicity of $G_{1},G_{2},\cdot \cdot \cdot, G_{N}$ provided that the number of nodes of each digraph $G_{i}$ is less than or equal to $3$ ($1\leq i\leq N$). The main result of paper is Theorem \ref{pende-bn}, where we completely solve the above question for tournaments. In addition to the above, we propose several open problems.

We use the standard notation henceforth. For any $s\in {\mathbb R},$ we set $\lfloor s \rfloor :=\sup \{
l\in {\mathbb Z} : s\geq l \}.$ Suppose that $X,\ Y,\ Z$ and $ T$ are given non-empty sets. Let us recall that a binary relation between $X$ into $Y$
is any subset
$\rho \subseteq X \times Y.$ 
If $\rho \subseteq X\times Y$ and $\sigma \subseteq Z\times T$ with $Y \cap Z \neq \emptyset,$ then
we define $\rho^{-1} \subseteq Y\times X$
and
$\sigma \circ \rho \subseteq X\times T$ by
$
\rho^{-1}:=\{ (y,x)\in Y\times X : (x,y) \in \rho \}
$
and
$$
\sigma \circ \rho :=\bigl\{(x,t) \in X\times T : \exists y\in Y \cap Z\mbox{ such that }(x,y)\in \rho\mbox{ and }
(y,t)\in \sigma \bigr\},
$$
respectively. Domain and range of $\rho$ are defined by $D(\rho):=\{x\in X :
\exists y\in Y\mbox{ such that }(x,y)\in \rho \}$ and $R(\rho):=\{y\in Y :
\exists x\in X\mbox{ such that }(x,y)\in \rho\},$ respectively; $\rho (x):=\{y\in Y : (x,y)\in \rho\}$ ($x\in X$), $ x\ \rho \ y \Leftrightarrow (x,y)\in \rho .$
Assuming $\rho$ is a binary relation on $X$ and $n\in {\mathbb N},$ we define $\rho^{n}
$ inductively; $\rho^{-n}:=(\rho^{n})^{-1}$ and $\rho^{0}:=
\{(x,x) : x\in X\}.$ Set $D_{\infty}(\rho):=\bigcap_{n\in {\mathbb N}} D(\rho^{n}),$ $\rho (X'):=\{y : y\in \rho(x)\mbox{ for some }x\in X'\}$ ($X'\subseteq X$) and
${\mathbb N}_{n}:=\{1,\cdot \cdot \cdot,n\}$ ($n\in {\mathbb N}$). By $P(A)$ and $\chi_{A}(\cdot),$ we denote the power set of $A$ and the characteristic function of $A,$ respectively. 

If $X$and $Y$ are topological spaces and $\rho \subseteq X\times Y,$ then we say that $\rho$ is continuous iff for every open subset $V$ of $Y$ there exists an open subset $U$ of $X$ such that $\rho^{-1}(V)=U \cap D(\rho).$ This clearly holds provided that $X$ is equipped with discrete topology.

\section{${\mathcal F}$-hypercyclicity and disjoint ${\mathcal F}$-hypercyclicity of binary relations: main definitions}\label{into}

Throughout the paper, we assume that $X$ and $Y$ are topological spaces as well as that $N\in {\mathbb N}$ and $N\geq 2.$ Suppose that ${\mathcal F}$ is a non-empty collection of certain subsets of ${\mathbb N},$ i.e., ${\mathcal F}\in P(P({\mathbb N}))$ and ${\mathcal F}\neq \emptyset.$ Observe that we do not require here that $|A|=\infty$ for all $A\in{\mathcal F}$ as well as that ${\mathcal F}$ satisfies the following property:
\begin{itemize}
\item[(I)] $B\in{\mathcal F}$ whenever there exists $A\in{\mathcal F}$ such that $A\subseteq B.$
\end{itemize}
Let us recall that, if ${\mathcal F}\in P(P({\mathbb N})) \setminus \emptyset $ satisfies (I), then it is said that 
${\mathcal F}$ is a Furstenberg family (\cite{furstenberg}); if so, then we say that ${\mathcal F}$ is a proper Furstenberg family iff $\emptyset \notin {\mathcal F}.$ For the sequel, we also need the notion of an upper Furstenberg family; that is any proper Furstenberg family ${\mathcal F}$ satisfying the following two conditions:
\begin{itemize}
\item[(II)] There exist a set $D$ and a countable set $M$ such that ${\mathcal F}=\bigcup_{\delta \in D} \bigcap_{\nu \in M}{\mathcal F}_{\delta,\nu},$ where for each $\delta \in D$ and $\nu \in M$ the following holds: If $A\in {\mathcal F}_{\delta,\nu},$ then there exists a
finite subset $F\subseteq {\mathbb N}$ such that the implication $A\cap F \subseteq B \Rightarrow B\in {\mathcal F}_{\delta,\nu}$ holds true.
\item[(III)] If $A\in {\mathcal F},$ then there exists $\delta \in D$ such that, for every $n\in {\mathbb N},$ we have $A-n\equiv \{k-n: k\in A,\ k>n\}\in {\mathcal F}_{\delta},$ where ${\mathcal F}_{\delta}\equiv \bigcap_{\nu \in M}{\mathcal F}_{\delta,\nu}.$
\end{itemize} 

We would like to propose the following definition (observe that the notion introduced can be further generalized by assuming that are given two non-empty families $\tau_{X} \in P(P(X))$ and $\tau_{Y}\in P(P(Y))$ satisfying that $(X,\tau_{X})$ and $(Y,\tau_{Y})$ are not necessarily topological spaces).

\begin{defn}\label{4-skins-MLO-okay}
Let $(\rho_{n})_{n\in {\mathbb N}}$ be a sequence of binary relations between the spaces
$X$ and $Y,$ let $\rho$ be a binary relation on $X$, and let $x\in X$. Suppose that ${\mathcal F}\in P(P({\mathbb N}))$ and ${\mathcal F}\neq \emptyset.$ Then we say
that:
\begin{itemize}
\item[(i)] $x$ is a strong ${\mathcal F}$-hypercyclic element of the sequence
$(\rho_{n})_{n\in {\mathbb N}}$ iff
$x\in \bigcap_{n\in {\mathbb N}} D(\rho_{n})$ and for each
$n\in {\mathbb N}$ there exists an element
$y_{n}\in \rho_{n}(x)$ such that for each open non-empty subset $V$ of $Y$ we have $\{ n\in {\mathbb N} : y_{n}\in V\}\in {\mathcal F} ;$ 
$(\rho_{n})_{n\in {\mathbb N}}$ is said to be strongly ${\mathcal F}$-hypercyclic iff there exists a strong ${\mathcal F}$-hypercyclic element of
$(\rho_{n})_{n\in {\mathbb N}}$;
\item[(ii)] $\rho$ is strong ${\mathcal F}$-hypercyclic iff the sequence
$(\rho^{n})_{n\in {\mathbb N}}$ is strong ${\mathcal F}$-hypercyclic; $x$ is said to be
a strong ${\mathcal F}$-hypercyclic element of $\rho$ iff $x$ is a strong ${\mathcal F}$-hypercyclic element of the sequence
$(\rho^{n})_{n\in {\mathbb N}};$
\item[(iii)] $x$ is an ${\mathcal F}$-hypercyclic element of the sequence
$(\rho_{n})_{n\in {\mathbb N}}$ iff
$x\in \bigcap_{n\in {\mathbb N}} D(\rho_{n})$ and for each open non-empty subset $V$ of $Y$ we have
$$
S(x,V):=\bigl\{ n\in {\mathbb N} : \rho_{n}x \cap V \neq \emptyset \bigr\}\in {\mathcal F} ;
$$ 
$(\rho_{n})_{n\in {\mathbb N}}$ is said to be ${\mathcal F}$-hypercyclic iff there exists an ${\mathcal F}$-hypercyclic element of
$(\rho_{n})_{n\in {\mathbb N}}$;
\item[(iv)] $\rho$ is ${\mathcal F}$-hypercyclic iff the sequence
$(\rho^{n})_{n\in {\mathbb N}}$ is ${\mathcal F}$-hypercyclic; $x$ is said to be
an ${\mathcal F}$-hypercyclic element of $\rho$ iff $x$ is an ${\mathcal F}$-hypercyclic element of the sequence
$(\rho^{n})_{n\in {\mathbb N}};$
\item[(v)] $(\rho_{n})_{n\in {\mathbb N}}$ is said to be strongly ${\mathcal F}$-topologically transitive iff for every open non-empty subset
$U\subseteq X$ and for every integer $n\in {\mathbb N}$ there exists an element $y_{n}\in \rho_{n}(U)$ such that for each open non-empty subset $V$ of $Y$ we have $\{ n\in {\mathbb N} : y_{n}\in V\}\in {\mathcal F} ;$ 
\item[(vi)] $\rho$ is strongly ${\mathcal F}$-topologically transitive iff the sequence $(\rho^{n})_{n\in {\mathbb N}}$ is strongly ${\mathcal F}$-topologically transitive;
\item[(vii)] $(\rho_{n})_{n\in {\mathbb N}}$ is said to be ${\mathcal F}$-topologically transitive iff for every two open non-empty subsets
$U\subseteq X$ and $V \subseteq Y$ we have 
$$
S(U,V):=\bigl\{ n\in {\mathbb N} : \rho_{n}(U)\cap V \neq \emptyset \bigr\}\in {\mathcal F} ;
$$ 
\item[(viii)] $\rho$ is ${\mathcal F}$-topologically transitive iff the sequence $(\rho^{n})_{n\in {\mathbb N}}$ is ${\mathcal F}$-topologically transitive.
\end{itemize}
\end{defn}

In any case set out above, the validity of (I) for ${\mathcal F}$ yields that the strong ${\mathcal F}$-hypercyclicity (topological transitivity) implies, in turn, the ${\mathcal F}$-hypercyclicity (topological transitivity) of considered sequence of binary relations (binary relation, element). 
This condition also ensures that, for every dynamical property introduced above, say ${\mathcal F}$-hypercyclicity, any extension of an ${\mathcal F}$-hypercyclic binary relation $\rho$ is likewise ${\mathcal F}$-hypercyclic (a similar statement holds for sequences of binary relations).  

The notion introduced in \cite{marina-prim} is recovered by setting that ${\mathcal F} $ is a collection of all non-empty subsets of ${\mathbb N}$ (in this case, generally, the notion of ${\mathcal F} $-hypercyclicity cannot be connected to that
of ${\mathcal F} $-topological transitivity in any reasonable way). It is worth noting 
that, in \cite{marina-prim}, the notion of strong ${\mathcal F}$-hypercyclicity and strong ${\mathcal F}$-topological transitivity (as well as their disjoint analogues) are called hypercyclicity and topological transitivity, respectively. So, the notions introduced in parts (iii)-(iv) of Definition \ref{4-skins-MLO-okay} as well as the notions introduced in parts (v)-(vi) of Definition \ref{4-skins-MLO-okay} are new.  

Definition \ref{4-skins-MLO-okay} is a rather general and covers some patological cases completely unambiguous to be further explored. Furthermore, the following holds:

\begin{itemize}
\item[(i)] The validity of (i), resp. (iii) [(ii), resp. (iv)], implies that $\bigcap_{n\in {\mathbb N}} D(\rho_{n}) \neq \emptyset$
[$D_{\infty}(\rho) \neq \emptyset$].
\item[(ii)] The validity of (v) [(vi)]
implies that  $D(\rho_{n}) \neq \emptyset$ for all $n\in {\mathbb N}$ [$D(\rho^{n}) \neq \emptyset$ for all $n\in {\mathbb N}$] but not $\bigcap_{n\in {\mathbb N}} D(\rho_{n}) \neq \emptyset$
[$D_{\infty}(\rho) \neq \emptyset$]. To illustrate this, consider first the case that $X=Y=\{x,y\}$ is equipped with discrete topology,
$\rho_{2n-1}:=\{(x,y)\},$ $\rho_{2n}:=\{(y,x)\}$ ($n\in {\mathbb N}$) and $\{{\mathbb N}, 2{\mathbb N},2{\mathbb N}+1\} \subseteq {\mathcal F}.$ Then it can be easily seen that the sequence $(\rho_{n})_{n\in {\mathbb N}}$ is strongly ${\mathcal F}$-topologically transitive,
$D(\rho_{n}) \neq \emptyset$ for all $n\in {\mathbb N}$ and $\bigcap_{n\in {\mathbb N}} D(\rho_{n}) = \emptyset .$ For (vi), it is sufficient to consider any binary relation $\rho $ on ${\mathbb N}=X=Y$ satisfying that $D(\rho^{n}) \neq \emptyset$ for all $n\in {\mathbb N}$ and $D_{\infty}(\rho) = \emptyset $; then   
we can endow $X$ and $Y$ with the anti-discrete topology $\tau=\{\emptyset,{\mathbb N}\}$ and $\rho$ will be strongly ${\mathcal F}$-topologically transitive provided that ${\mathbb N} \in {\mathcal F}.$
\item[(iii)]
In the case of consideration parts
(vii) and (viii), we do not need to have that $D(\rho_{n}) = \emptyset$ for all $n\in {\mathbb N}$ [$D(\rho^{n}) \neq \emptyset $ for all $n\in {\mathbb N}$]; to see this, assume that $X=Y=\{x,y\}$ is equipped with discrete topology, $\rho :=\{(x,y)\}$ and $\{\emptyset , \{1\}\}\subseteq {\mathcal F}.$ Then $\rho$ is (not strongly) ${\mathcal F}$-topologically transitive, $D(\rho^{n}) = \emptyset$ for all $n\geq 2$ and therefore
$D_{\infty}(\rho)=\emptyset .$ 
\item[(iv)] Assume ${\mathbb N} \notin {\mathcal F}$ and $D(\rho_{n})\neq \emptyset$ for all $n\in {\mathbb N}$ [$D_{\infty}(\rho) \neq \emptyset$]. Then we can easily seen by plugging ($U=X$ and) $V=Y$ that the sequence $(\rho_{n})_{n\in {\mathbb N}}$ [$\rho$] cannot satisfy any of the above introduced properties. 
\end{itemize}

If $X=Y$ and $(\rho_{n})_{n\in {\mathbb N}}$ is a sequence of symmetric binary relations on $X,$ then for each $x\in \bigcap_{n\in {\mathbb N}} D(\rho_{n})$ and for each open non-empty subset $V$ of $Y$ we have 
$ S(x,V) +2{\mathbb N}\subseteq S(x,V) ,$ so that $(\rho_{n})_{n\in {\mathbb N}}$ cannot be ${\mathcal F}$-hypercyclic if for each subset $A\in {\mathcal F}$ the assumption $A\neq \emptyset$ implies that $A+2{\mathbb N}$ is not a subset of $A;$ furthermore, in this case, for every two open non-empty subsets
$U\subseteq X$ and $V \subseteq X,$ we have
$
S(U,V)=S(V,U).$

\begin{rem}\label{qselo}
Assume that $x$ is a (strong) ${\mathcal F}$-hypercyclic element of a binary relation $\rho$ on $X,$ $l\in {\mathbb N}$ and $x\in \rho^{l}z$ for some element $z\in X.$ If (for each open non-empty subset $V$ of $Y$ and for each sequence $(\omega_{n})$ in $X$) for each open non-empty subset $V$ of $Y$ and for each $\omega\in \bigcap_{n\in {\mathbb N}} D(\rho_{n})$ the supposition ($\{n\in {\mathbb N} : \omega_{n+l}\in V\}\in {\mathcal F}$ implies $\{n\in {\mathbb N} : \omega_{n}\in V\} \in {\mathcal F}$) $\{n\in {\mathbb N} : \rho^{n+l}\omega\in V\}\in {\mathcal F}$ implies $\{n\in {\mathbb N} : \rho^{n}\omega\in V\} \in {\mathcal F},$ then $z$ is likewise a (strong) ${\mathcal F}$-hypercyclic element for $\rho .$ This, in particular, holds if ${\mathcal F}$ is a collection of all subsets of ${\mathbb N}$ which do have at least $m$ elements, where $m\in {\mathbb N}_{0}.$ 
\end{rem}

In the following two definitions, we consider disjoint analogues of notions introduced in Definition \ref{4-skins-MLO-okay}:

\begin{defn}\label{4-skins-MLO-okay-novio}
Suppose that ${\mathcal F}\in P(P({\mathbb N})),$ ${\mathcal F}\neq \emptyset,$ $N\geq 2,$ $(\rho_{j,n})_{n\in {\mathbb N}}$
is a sequence of binary relations between the spaces $X$ and $Y$ ($1\leq j\leq N$),
$\rho_{j}$ is a binary relation on $X$ ($1\leq j\leq N$) and $x\in X$. Then we say
that:
\begin{itemize}
\item[(i)] $x$ is a strong $d{\mathcal F}$-hypercyclic element of the sequences
$(\rho_{1,n})_{n\in {\mathbb N}},\cdot \cdot \cdot , (\rho_{N,n})_{n\in {\mathbb N}}$
iff for each $n\in {\mathbb N}$ there exist elements $y_{j,n}\in \rho_{j,n}(x)$ ($1\leq j\leq N$) such that for every open non-empty subsets $V_{1},\cdot\cdot \cdot , V_{N}$ of $Y,$ we have $\{n\in {\mathbb N} : y_{1,n}\in V_{1},\ y_{2,n}\in V_{2},\cdot \cdot \cdot,\ y_{N,n}\in V_{N}\}\in {\mathcal F};$ 
the sequences $(\rho_{1,n})_{n\in {\mathbb N}},\cdot \cdot \cdot ,
(\rho_{N,n})_{n\in {\mathbb N}}$ are called
strongly $d{\mathcal F}$-hypercyclic iff there exists a strong $d{\mathcal F}$-hypercyclic element of
$(\rho_{1,n})_{n\in {\mathbb N}},\cdot \cdot \cdot , (\rho_{N,n})_{n\in {\mathbb N}};$
\item[(ii)] $x$ is a strong $d{\mathcal F}$-hypercyclic element of the binary relations $\rho_{1},\cdot \cdot \cdot ,
\rho_{N}$ iff $x$ is a strong $d{\mathcal F}$-hypercyclic element of the sequences
$(\rho_{1}^{n})_{n\in {\mathbb N}},\cdot \cdot \cdot , (\rho_{N}^{n})_{n\in {\mathbb N}};$
the binary relations $\rho_{1},\cdot \cdot \cdot , \rho_{N}$
are called strongly $d{\mathcal F}$-hypercyclic iff there exists a strong $d{\mathcal F}$-hypercyclic element of
$\rho_{1},\cdot \cdot \cdot , \rho_{N};$
\item[(iii)] $x$ is a $d{\mathcal F}$-hypercyclic element of the sequences
$(\rho_{1,n})_{n\in {\mathbb N}},\cdot \cdot \cdot , (\rho_{N,n})_{n\in {\mathbb N}}$
iff $x\in \bigcap_{1\leq j \leq N, n\in {\mathbb N}} D_{\infty}(\rho_{j,n})$ and for every open non-empty subsets $V_{1},\cdot\cdot \cdot , V_{N}$ of $Y,$ we have (${\mathrm V}=(V_{1},\ V_{2},\cdot \cdot \cdot, V_{N})$)
$$
S(x,{\mathrm V}):=\bigl\{ n\in {\mathbb N} : \rho_{1,n}x \cap V_{1}\neq \emptyset, \ \rho_{2,n}x \cap V_{2}\neq \emptyset,  \cdot \cdot \cdot ,\ \rho_{N,n}x \cap V_{N} \neq \emptyset \bigr\}\in {\mathcal F} ;
$$ 
the sequences $(\rho_{1,n})_{n\in {\mathbb N}},\cdot \cdot \cdot ,
(\rho_{N,n})_{n\in {\mathbb N}}$ are called
$d{\mathcal F}$-hypercyclic iff there exists a $d{\mathcal F}$-hypercyclic element of
$(\rho_{1,n})_{n\in {\mathbb N}},\cdot \cdot \cdot , (\rho_{N,n})_{n\in {\mathbb N}};$
\item[(iv)] $x$ is a $d{\mathcal F}$-hypercyclic element of the binary relations $\rho_{1},\cdot \cdot \cdot ,
\rho_{N}$ iff $x$ is a $d{\mathcal F}$-hypercyclic element of the sequences
$(\rho_{1}^{n})_{n\in {\mathbb N}},\cdot \cdot \cdot , (\rho_{N}^{n})_{n\in {\mathbb N}};$
the binary relations $\rho_{1},\cdot \cdot \cdot , \rho_{N}$
are called $d{\mathcal F}$-hypercyclic iff there exists a $d{\mathcal F}$-hypercyclic element of
$\rho_{1},\cdot \cdot \cdot , \rho_{N}.$
\end{itemize}
\end{defn}

\begin{defn}\label{4-skins-MLO-okay-novio'frim}
Suppose that ${\mathcal F}\in P(P({\mathbb N})),$ ${\mathcal F}\neq \emptyset,$ $N\geq 2,$ $(\rho_{j,n})_{n\in {\mathbb N}}$
is a sequence of binary relations between the spaces $X$ and $Y$ ($1\leq j\leq N$), and
$\rho_{j}$ is a binary relation on $X$  ($1\leq j\leq N$). Then we say
that:
\begin{itemize}
\item[(i)] the sequences $(\rho_{1,n})_{n\in {\mathbb N}},\cdot \cdot \cdot , (\rho_{N,n})_{n\in {\mathbb N}}$
are strongly $d{\mathcal F}$-topologically transitive iff 
for every open non-empty subset
$U\subseteq X$ and 
for every open non-empty subsets $V_{1},\cdot\cdot \cdot , V_{N}$ of $Y,$ 
there exists an element $x\in U$ such that, 
for every integers $n\in {\mathbb N}$ and $j\in {\mathbb N}_{N},$ there exists an element $y_{j,n}\in \rho_{j,n}x$ so that $\{ n\in {\mathbb N} : y_{j,n}\in V_{j}\mbox{ for all }j\in {\mathbb N}_{N}\}\in {\mathcal F} ;$ 
\item[(ii)] the binary relations $\rho_{1},\cdot \cdot \cdot , \rho_{N}$ are called
strongly $d{\mathcal F}$-topologically transitive iff the sequences
$(\rho_{1}^{n})_{n\in {\mathbb N}},\cdot \cdot \cdot , (\rho_{N}^{n})_{n\in {\mathbb N}}$
are strongly $d{\mathcal F}$-topologically transitive;
\item[(iii)] the sequences $(\rho_{1,n})_{n\in {\mathbb N}},\cdot \cdot \cdot , (\rho_{N,n})_{n\in {\mathbb N}}$
are $d{\mathcal F}$-topologically transitive iff 
for every open non-empty subset
$U\subseteq X$ and 
for every open non-empty subsets $V_{1},\cdot\cdot \cdot , V_{N}$ of $Y,$ 
we have 
$\{ n\in {\mathbb N} : (\exists x \in U)\, \rho_{j,n}x \cap V_{j} \neq \emptyset\mbox{ for all }j\in {\mathbb N}_{N}\}\in {\mathcal F} ;$ 
\item[(iv)] the binary relations $\rho_{1},\cdot \cdot \cdot , \rho_{N}$ are
$d{\mathcal F}$-topologically transitive iff the sequences
$(\rho_{1}^{n})_{n\in {\mathbb N}},\cdot \cdot \cdot , (\rho_{N}^{n})_{n\in {\mathbb N}}$
are $d{\mathcal F}$-topologically transitive.
\end{itemize}
\end{defn}

If
the binary relations $\rho_{1},\cdot \cdot \cdot , \rho_{N}$ are (strongly) $d{\mathcal F}$-hypercyclic 
((strongly) $d{\mathcal F}$-topologically transitive), then we also say that the tuple $(\rho_{1},\cdot \cdot \cdot , \rho_{N})$ is strong $d{\mathcal F}$-hypercyclic 
((strongly) $d{\mathcal F}$-topologically transitive) and vice versa.

We have the following simple observations:

\begin{itemize}
\item[(i)] Definition \ref{4-skins-MLO-okay-novio}:  The validity of (i), resp. (iii) [(ii), resp. (iv)], implies that $\bigcap_{n\in {\mathbb N},j\in {\mathbb N}_{N}} D(\rho_{j,n}) \neq \emptyset$
[$\bigcap_{j\in {\mathbb N}_{N}} D_{\infty}(\rho_{j}) \neq \emptyset$].
\item[(ii)]  Definition \ref{4-skins-MLO-okay-novio'frim}: The validity of (i) [(ii)]
implies that $\bigcap_{j\in {\mathbb N}_{N}}D(\rho_{j,n}) \neq \emptyset$ for all $n\in {\mathbb N}$ [$\bigcap_{j\in {\mathbb N}_{N}} D(\rho_{j}^{n}) \neq \emptyset$ for all $n\in {\mathbb N}$] but not $\bigcap_{n\in {\mathbb N}} D(\rho_{j,n}) \neq \emptyset$ for some $j\in {\mathbb N}_{N}$
[$D_{\infty}(\rho_{j}) \neq \emptyset$ for some $j\in {\mathbb N}_{N}$]. 
\item[(iii)] Definition \ref{4-skins-MLO-okay-novio'frim}: The validity of (iii) [(iv)]
does not imply that there exists $j\in {\mathbb N}_{N}$ such that $D(\rho_{j,n}) \neq \emptyset$ for all $n\in {\mathbb N}$ [$D(\rho_{j}^{n}) \neq \emptyset $ for all $n\in {\mathbb N}$]; to verify this, let $X=\{x,y\}$ be equipped with topology $\tau_{1}=\{\emptyset,\{x\},\{x,y\}\},$ let $Y=\{x,y\}$ be equipped with topology $\tau_{2}=\{\emptyset,\{y\},\{x,y\}\},$ and let $\rho_{1}=\rho_{2}=\{(x,y)\}.$ Suppose that
${\mathcal F}=P(P({\mathbb N})) \setminus \{\emptyset\}.$
Then $D(\rho_{1}^{n})=\emptyset$ for $n\geq 2,$ $\rho_{1}$ and $\rho_{2}$ are $d{\mathcal F}$-topologically transitive but not strongly $d{\mathcal F}$-topologically transitive (see also Remark \ref{prcojed} below).
\item[(iv)] Assume ${\mathbb N} \notin {\mathcal F}.$ If $\bigcap_{j\in {\mathbb N}_{N}} D(\rho_{j,n})\neq \emptyset$ for all $n\in {\mathbb N}$ [$\bigcap_{j\in {\mathbb N}_{N}}D(\rho_{j}^{n}) \neq \emptyset$ for all $n\in {\mathbb N}$], then the sequences $(\rho_{1,n})_{n\in {\mathbb N}},\cdot \cdot \cdot , (\rho_{N,n})_{n\in {\mathbb N}}$ [binary relations $\rho_{1},\cdot \cdot \cdot , \rho_{N}$] cannot satisfy any of the above introduced disjoint properties. 
\end{itemize}

\begin{rem}\label{topologije-prc}
\begin{itemize}
\item[(i)] In the parts (i) and (iii) of Definition \ref{4-skins-MLO-okay-novio}, the topology on $X$ does not play any role. We can assume that $Y$ is equipped with arbitrary topologies $\tau_{1}^{Y},\cdot \cdot \cdot,\tau_{N}^{Y}$ and that, for every $i\in {\mathbb N}_{N},$ $V_{i}$ is open for the topology $\tau_{i}^{Y}.$    
\item[(ii)] In the parts (ii) and (iv) of Definition \ref{4-skins-MLO-okay-novio}, we can assume that $X=Y$ is equipped with arbitrary topologies $\tau_{1}^{X},\cdot \cdot \cdot,\tau_{N}^{X}$ and that, for every $i\in {\mathbb N}_{N},$ $V_{i}$ is open for the topology $\tau_{i}^{X}.$
\end{itemize}
Similar observations can be given for Definition \ref{4-skins-MLO-okay-novio'frim}; albeit a great number of our results continues to hold with this extended notion, we will analyze henceforth only the usually considered case that $X$ is equipped with exactly one topology and $Y$ is  equipped with exactly one topology.
\end{rem}

We round off this section by stating the following simple proposition, stated here without a corresponding proof which can be left to the interested readers:

\begin{prop}\label{komponente}
\begin{itemize}
\item[(i)]
Suppose that the sequences $(\rho_{1,n})_{n\in {\mathbb N}},\cdot \cdot \cdot , (\rho_{N,n})_{n\in {\mathbb N}}$ of binary relations between the spaces $X$ and $Y$, resp. the
binary relations $\rho_{1},\rho_{2}, \cdot \cdot \cdot , \rho_{N}$ on $X,$ are $d{\mathcal F}$-hypercyclic ($d{\mathcal F}$-topologically transitive). Then for each $j\in {\mathbb N}_{N}$ 
the sequence $(\rho_{j,n})_{n\in {\mathbb N}},$ resp. the binary relation $\rho_{j},$
is ${\mathcal F}$-hypercyclic (${\mathcal F}$-topologically transitive) provided that $(\rho_{1,n})_{n\in {\mathbb N}}=(\rho_{2,n})_{n\in {\mathbb N}}=\cdot \cdot \cdot = (\rho_{N,n})_{n\in {\mathbb N}}$, resp. $\rho_{1}=\rho_{2}=\cdot \cdot \cdot = \rho_{N},$ or that the condition \emph{(I)} holds for ${\mathcal F}.$
\item[(ii)] Suppose that the sequences $(\rho_{1,n})_{n\in {\mathbb N}},\cdot \cdot \cdot , (\rho_{N,n})_{n\in {\mathbb N}}$ of binary relations between the spaces $X$ and $Y$, resp. the
binary relations $\rho_{1},\rho_{2}, \cdot \cdot \cdot , \rho_{N}$ on $X,$ are strongly $d{\mathcal F}$-hypercyclic (strongly $d{\mathcal F}$-topologically transitive). Then for each $j\in {\mathbb N}_{N}$ 
the sequence $(\rho_{j,n})_{n\in {\mathbb N}},$ resp. the binary relation $\rho_{j},$
is strongly ${\mathcal F}$-hypercyclic (strongly ${\mathcal F}$-topologically transitive) provided that the condition \emph{(I)} holds for ${\mathcal F}.$
\end{itemize}
\end{prop}

\section{Results for general binary relations, simple graphs and digraphs}\label{guptal}

As already mentioned, we will divide this section into three separate subsections. In the first one, we will present a few result about ${\mathcal F}$-hypercyclicity of general binary relations.

\subsection{Results for general binary relations}\label{general binary} 

For the sequel, set ${\mathrm D}:=\bigcap_{n\in {\mathbb N}} D(\rho_{n})$ and $\check{\rho_{n}}:=\{(x,y)\in \rho_{n} : x\in {\mathrm D}\},$ $n\in {\mathbb N}$ [$\check{\rho}:=\{(x,y)\in \rho : x\in D_{\infty}(\rho)\}$]. Then $\check{\rho_{n}}$ [$\check{\rho}$] is a binary relation between ${\mathrm D}$ and $Y,$ with $D(\check{\rho_{n}})={\mathrm D}$ for all $n\in {\mathbb N}$ [$D_{\infty}(\rho)$ and $Y,$ with $D(\check{\rho})=D_{\infty}(\rho)$]. 

The following proposition holds true: 

\begin{prop}\label{az}
Let $(\rho_{n})_{n\in {\mathbb N}}$ be a sequence of binary relations between the spaces
$X$ and $Y,$ let $\rho$ be a binary relation on $X$, and let $x\in X$. Then $x$ is a (strong) ${\mathcal F}$-hypercyclic element of the sequence
$(\rho_{n})_{n\in {\mathbb N}}$ iff $x$ is a (strong) ${\mathcal F}$-hypercyclic element of the sequence
$(\check{\rho_{n}})_{n\in {\mathbb N}};$ in particular, $(\rho_{n})_{n\in {\mathbb N}}$ is (strongly) ${\mathcal F}$-hypercyclic iff $(\check{\rho_{n}})_{n\in {\mathbb N}}$ is (strongly) ${\mathcal F}$-hypercyclic.
\end{prop} 

In a certain sense, the above proposition shows that it is sufficient to introduce the notions of (strong) ${\mathcal F}$-hypercyclicity only for binary relations whose domain is the whole space $X.$ But, this is actually not the case because we need to know some further properties of ${\mathrm D}$ in $X;$ for example, in the next generalization
of implications (a) $\Rightarrow$ (b) $\Rightarrow$ (c) $\Rightarrow$ (d) of \cite[Theorem 15]{boni-upper}, we impose the condition that the subspace ${\mathrm D}$ of $X$ is a Baire space, which particularly holds in the following two special cases: $X$ is a Baire space and ${\mathrm D}$ is open in $X$ or $X$ is a complete metric space and ${\mathrm D}$ is a closed subspace of $X:$

\begin{thm}\label{isterace}
Let $(\rho_{n})_{n\in {\mathbb N}}$ be a sequence of binary relations between the topological spaces
$X$ and $Y,$ let the subspace ${\mathrm D}$ of $X$ be a Baire space, and let $Y$ be second-countable. Assume that $\check{\rho_{n}}\subseteq {\mathrm D} \times Y$ is continuous for all $n\in {\mathbb N}.$ If ${\mathcal F}$ is a Furstenberg family and \emph{(II)} holds, then we have \emph{(i)} $\Rightarrow$ \emph{(ii)} $\Rightarrow$ \emph{(iii)} $\Rightarrow$ \emph{(iv)}, where:
\begin{itemize}
\item[(i)] For any non-empty open subset $V$ of $Y$ there is some $\delta \in D$ such that for any non-empty open subset $U$ of $X$ such that $U\cap {\mathrm D}\neq \emptyset$ there is some $x\in U\cap {\mathrm D}$ such that
$
\{ n\in {\mathbb N} : \rho_{n}x \cap V \neq \emptyset \} \in {\mathcal F}_{\delta}.
$
\item[(ii)] For any non-empty open subset $V$ of $Y$ there is some $\delta \in D$ such that, for any non-empty open subset $U$ of $X$ such that $U\cap {\mathrm D}\neq \emptyset$ and for every $\nu \in M$ there is some $x\in U\cap {\mathrm D}$ such that
$
\{ n\in {\mathbb N} : \rho_{n}x \cap V \neq \emptyset \} \in {\mathcal F}_{\delta , \nu}.
$
\item[(iii)] The set consisting of all ${\mathcal F}$-hypercyclic vectors of $(\rho_{n})_{n\in {\mathbb N}}$ is residual in ${\mathrm D}.$
\item[(iv)] The sequence $(\rho_{n})_{n\in {\mathbb N}}$ is ${\mathcal F}$-hypercyclic.
\end{itemize}
\end{thm}

\proof
The implications (i) $\Rightarrow$ (ii) and (iii) $\Rightarrow$ (iv) are trivial and all that we need to show is that (ii) implies (iii). For this, we can repeat almost literally the arguments given in the proof of corresponding implication (b) $\Rightarrow$ (c) of \cite[Theorem 15]{boni-upper}, with the sequence $(T_{n})$ and term $T_{n}x\in V_{k}$ replaced therein with the sequence of continuous relations 
$(\check{\rho_{n}})_{n\in {\mathbb N}}$ and term $\rho_{n}x \cap V_{k} \neq \emptyset,$ showing that the set consisting of all ${\mathcal F}$-hypercyclic vectors of $(\check{\rho_{n}})_{n\in {\mathbb N}}$ is residual in ${\mathrm D}.$ After that, we can apply Proposition \ref{az}.
\endproof

In contrast to \cite{boni-upper}, we do not use the conditions that $X$ and $Y$ are metric spaces, as well as the condition (III). 
Keeping in mind  Proposition \ref{az} and Theorem \ref{isterace}, it seems plausible that the assertions of ${\mathcal A}-$Hypercyclicity Criterion \cite[Theorem 22]{boni-upper} and ${\mathcal A}-$Hypercyclicity Criterion, second version \cite[Theorem 26]{boni-upper}, can be extended for continuous multivalued linear operators (see \cite{marina} for the notion of a multivalued linear operator; the continuity is understood in the sense of continuity of a general binary relation). It also seems plausible that a great number of other Hypercyclicity Criteria known in the existing literature can be formulated for continuous multivalued linear operators. We will not discuss these questions in more detail here.

Now we will turn our attention in another direction, by giving a few useful observations in the case that $X$ is equipped with discrete topology or anti-discrete topology. Suppose first that $X$ carries the anti-discrete topology $\tau=\{\emptyset , X\}$ and binary relations $\rho,$ $\rho_{1},\cdot \cdot \cdot, \rho_{N}$ on $X$ are given. Then any vector $x\in D_{\infty}(\rho)$ is a strong ${\mathcal F}$-hypercyclic vector for $\rho$ (strong $d{\mathcal F}$-hypercyclic vector for $\rho_{1},\cdot \cdot \cdot, \rho_{N}$), so that the notions of ${\mathcal F}$-hypercyclicity and strong ${\mathcal F}$-hypercyclicity ($d{\mathcal F}$-hypercyclicity and strong $d{\mathcal F}$-hypercyclicity) coincide; the same holds for the notions of ${\mathcal F}$-topological transitivity and strong ${\mathcal F}$-topological transitivity ($d{\mathcal F}$-topological transitivity and strong $d{\mathcal F}$-topological transitivity). This is no longer true in the case that $\tau$ is not the anti-discrete topology and we will illustrate this only for the ${\mathcal F}$-hypercyclicity: let $X=\{x_{1},x_{2}\},$ $\rho=\{(x_{1},x_{2}),(x_{2},x_{2})\}$ and $\tau =\{\emptyset, \{x_{2}\},\{x_{1},x_{2}\}\}.$ Then $x_{1}$ and $x_{2}$ are both ${\mathcal F}$-hypercyclic vectors for $\rho,$ while $x_{2}$ is the only strong ${\mathcal F}$-hypercyclic
vector for $\rho.$ 

If $X=\{x_{1},x_{2},\cdot \cdot \cdot, x_{n}\}$ carries the discrete topology, then
an element
$x\in X$ is an ${\mathcal F}$-hypercyclic vector for a binary relation $\rho$ on $X$ iff for every non-empty subset $V$ of ${\mathbb N}_{n}$ we have $\{k\in {\mathbb N} : (\exists i\in V ) \, x_{i}\in \rho^{k}x\} \in {\mathcal F}.$ Since
$$
\bigl\{k\in {\mathbb N} : (\exists i\in V ) \, x_{i}\in \rho^{k}x\bigr\} =\bigcup_{i\in V} \bigl\{k\in {\mathbb N} : x_{i}\in \rho^{k}x\bigr\},
$$
we have the following: Assume that ${\mathcal F}$ is closed under finite unions. Then $x\in X$ is an ${\mathcal F}$-hypercyclic vector for $\rho$ iff $x\in D_{\infty}(\rho)$ and for each $i\in {\mathbb N}_{n}$ we have $\{k\in {\mathbb N} : x_{i}\in \rho^{k}x\}\in {\mathcal F}.$ Arguing similarly we can prove that, under the same assumption on ${\mathcal F},$ $x\in X$ is a $d{\mathcal F}$-hypercyclic vector for binary relations $\rho_{1},\rho_{2},\cdot \cdot \cdot,\rho_{N}$ on $X$ iff $x\in D_{\infty}(\rho_{j})$ for $1\leq j\leq N$ and for 
any choice of elements $x_{i_{1}},x_{i_{2}},\cdot \cdot \cdot,x_{i_{N}}$ in $X$ ($1\leq i_{s}\leq n$ for $s\in {\mathbb N}_{N}$) we have $\bigcap_{j\in {\mathbb N}_{N}} \{k\in {\mathbb N} : x_{i_{j}}\in \rho_{j}^{k}x\}\in {\mathcal F}.$

Consider the case that $X=Y=\{x_{1},x_{2},\cdot \cdot \cdot, x_{n}\}$ is equipped with arbitrary topology. For any binary relation $\rho$ on $X,$
by $[\rho]$ we denote the adjacency matrix of $\rho,$ defined by $a_{ij}:=1$ if $x_{i} \, \rho\,  x_{j}$ and $a_{ij}:=0,$ otherwise. By a $\rho$-walk, we mean any finite sequence $x_{i_{1}}x_{i_{2}}\cdot \cdot \cdot x_{i_{s}},$ where $s\in {\mathbb N} \setminus \{1\},$ $1\leq i_{j} \leq n$ for $1\leq j\leq s$ and $x_{i_{j}} \, \rho \, x_{i_{j+1}}$ for $1\leq j\leq s-1;$ the length of $x_{i_{1}}x_{i_{2}}\cdot \cdot \cdot x_{i_{s}}$ is said to be $s$, while $x_{i_{1}}$ and $x_{i_{s}}$ are said to be the starting and ending point of $x_{i_{1}}x_{i_{2}}\cdot \cdot \cdot x_{i_{s}}.$  We also say that $x_{i_{1}}x_{i_{2}}\cdot \cdot \cdot x_{i_{s}}$ is an $(x_{i_{1}}-x_{i_{s}})$ $\rho$-walk. 
Set $[\rho]^{k}:=[a_{i,j}^{k}]_{1\leq i,j\leq n}$ ($k\in {\mathbb N}$). 
Arguing as in the case of simple graphs (\cite{petrovic}), we can simply conclude that the number of different $(x_{i}-x_{j})$ $\rho$-walks of length $k$ equals $a_{i,j}^{k}$ ($1\leq i,j\leq n,$ $k\in {\mathbb N}$). This fact enables one to simply reformulate the notion introduced in Definition \ref{4-skins-MLO-okay} in terms of appropriate conditions involving the adjacency matrix $[\rho]:$ 

\begin{prop}\label{boze}
\begin{itemize}
\item[(i)] $x_{i}$ is an ${\mathcal F}$-hypercyclic vector for $\rho$ ($i\in {\mathbb N}_{n}$) iff 
for every $k\in {\mathbb N}$ there exists a $\rho$-walk of length $k$ starting at $x_{i}$
and for each open non-empty subset $V$ of $X$
we have $\{ k\in {\mathbb N} : (\exists j\in V)\, a_{ij}^{k}\geq 1\} \in {\mathcal F}.$ 
\item[(ii)] $\rho$ is ${\mathcal F}$-topologically transitive iff for each pair of open non-empty subsets $U,\ V$ of $X$ we have $\{ k\in {\mathbb N} : (\exists i\in U)\, (\exists j\in V)\, a_{ij}^{k}\geq 1\} \in {\mathcal F}.$
\end{itemize}
\end{prop}

The situation is quite similar for disjointness. If $[\rho_{s}]=[a_{ij}^{s}]_{1\leq i,j\leq n}$ is the adjacency matrix of a binary relation $\rho_{s}$ on $X$, then we denote $[\rho_{s}]^{k}=[a_{ij}^{k;s}]_{1\leq i,j\leq n}$ ($s\in {\mathbb N}_{N},$ $k\in {\mathbb N}$). We have the following:  

\begin{prop}\label{boze-prim}
Let $\rho_{s}$ be a binary relation on $X$ ($1\leq s\leq N$).
\begin{itemize}
\item[(i)] $x_{i}$ is a $d{\mathcal F}$-hypercyclic vector for $\rho_{1}, \rho_{2},\cdot \cdot \cdot, \rho_{N}$ ($i\in {\mathbb N}_{n}$) iff 
for every $k\in {\mathbb N}$ and $s\in {\mathbb N}_{N}$ there exists a $\rho_{s}$-walk of length $k$ starting at $x_{i}$
and for each open non-empty subsets $V_{s}$ of $ X_{s}$ ($1\leq s\leq N$)
we have $\{ k\in {\mathbb N} : (\forall s\in {\mathbb N}_{N})(\exists j_{s}\in V_{s})\, a_{ij_{s}}^{k;s}\geq 1\} \in {\mathcal F}.$ 
\item[(ii)] $\rho_{1}, \rho_{2},\cdot \cdot \cdot, \rho_{N}$ are $d{\mathcal F}$-topologically transitive iff  for each open non-empty subsets $U,V_{1},V_{2},\cdot \cdot \cdot,V_{N}$ of $ X$ we have $\{ k\in {\mathbb N} : (\exists i\in U)(\forall s\in {\mathbb N}_{N})(\exists j_{s}\in V_{s}) a_{ij_{s}}^{k;s}\geq 1\} \in {\mathcal F}.$
\end{itemize}
\end{prop}

\subsection{Results for simple graphs}\label{graphs} 

Let $X=G=\{x_{1},x_{2},\cdot \cdot \cdot, x_{n}\}$ be finite, let $|G|>1,$ and let $\rho$ be a symmetric relation on $G$ such that, for every $g\in G,$ we have
$(g,g)\notin \rho .$ As it is well-known, $(G,\rho)$ is said to be a simple graph (see \cite{bondy}, \cite{booksa} and \cite{petrovic} for the basic theory of graphs). By $E(G)$ we denote the set consisting of all unoriented arcs of $G.$ The notion of distance $d(u,v)$ of two nodes $u,\ v\in G,$ as well as the notions of diameter $d(G)$ of graph $G,$ walks, paths and closed contours in $G$ are defined usually (let us only recall that $d(u,u)=0,$ $u\in G$). By $[A(G)]$ we denote the adjacency matrix of $G.$ For more details about applications of matrix theory to graphs, we refer the reader to the monographs \cite{dragos}-\cite{dragos1}.

Suppose, for the time being, that ${\mathcal F}=P(P({\mathbb N})) \setminus \{\emptyset \}.$ Then it can be easily seen that the graph $G,$ equipped with discrete topology, is connected iff $G$ is (strongly) ${\mathcal F}$-hypercyclic iff $G$ is (strongly) ${\mathcal F}$-topologically transitive (\cite{marina-prim}); if this is the case, then any element of $G$ is a (strong) hypercyclic element of $\rho.$ 
Furthermore, if $G$ is equipped with arbitrary topology, then $G$ is ${\mathcal F}$-hypercyclic (${\mathcal F}$-topologically transitive) iff $G$ is strongly ${\mathcal F}$-hypercyclic  (strongly ${\mathcal F}$-topologically transitive).
For disjointness, a similar statement holds true:

\begin{prop}\label{poka}
Let ${\mathcal F}=P(P({\mathbb N})) \setminus \{\emptyset \}.$
Suppose that $G_{1},G_{2},\cdot \cdot \cdot, G_{N}$ are given graphs with the same set of nodes $X=\{x_{1},x_{2},\cdot \cdot \cdot, x_{n}\}$. Then $G_{1},G_{2},\cdot \cdot \cdot, G_{N}$ are $d{\mathcal F}$-hypercyclic ($d{\mathcal F}$-topologically transitive) iff $G_{1},G_{2},\cdot \cdot \cdot, G_{N}$ are strongly $d{\mathcal F}$-hypercyclic (strongly $d{\mathcal F}$-topologically transitive).
\end{prop}

\proof 
We will prove the statement only for $d{\mathcal F}$-hypercyclicity and strong $d{\mathcal F}$-hypercyclicity.
Since the condition (I) holds, we only need to show that $d{\mathcal F}$-hypercyclicity of $G_{1},G_{2},\cdot \cdot \cdot, G_{N}$ implies their strong $d{\mathcal F}$-hypercyclicity. Let $x\in X$ be a $d{\mathcal F}$-hypercyclic vector of $G_{1},G_{2},\cdot \cdot \cdot, G_{N};$ we will prove that $x\in X$ is a strong $d{\mathcal F}$-hypercyclic vector of $G_{1},G_{2},\cdot \cdot \cdot, G_{N}.$ Let $a_{n}$ denote the number of open non-empty subsets of $G_{i},$ which will be denoted by $V_{s}$ ($1\leq s \leq a_{n}$). Set $r_{n}:=a_{n}^{N}.$ Let the tuples $(V_{1},V_{1},\cdot \cdot \cdot,V_{1}),\cdot \cdot \cdot, (V_{a_{n}},V_{a_{n}},\cdot \cdot \cdot,V_{a_{n}})$
be listed in some alphabetic order.
For the tuple $(V_{1},V_{1},\cdot \cdot \cdot,V_{1}),$ we know that there exist a positive integer $k_{1}\in {\mathbb N}$ and elements $y_{1,k_{1}}\in \rho_{1}^{k_{1}}x \cap V_{1},\cdot \cdot \cdot , y_{N,k_{1}}\in \rho_{N}^{k_{1}}x \cap V_{1}.$ By the symmetry of relations $\rho_{1},\rho_{2},\cdot \cdot \cdot, \rho_{N},$ for the tuple $(V_{1},V_{1},\cdot \cdot \cdot,V_{1}, V_{2}),$ there exist a positive integer $k_{2}>k_{1}$ and elements $y_{1,k_{2}}\in \rho_{1}^{k_{2}}x \cap V_{1},\cdot \cdot \cdot ,y_{N-1,k_{2}}\in \rho_{N-1}^{k_{2}}x \cap V_{1},  y_{N,k_{2}}\in \rho_{N}^{k_{2}}x \cap V_{2}.$ Repeating this procedure, for the tuple $(V_{a_{n}},V_{a_{n}},\cdot \cdot \cdot,V_{a_{n}}),$ there exist a positive integer $k_{r_{n}}>k_{r_{n}-1}$ and elements $y_{1,k_{r_{n}}}\in \rho_{1}^{k_{r_{n}}}x \cap V_{a_{n}},\cdot \cdot \cdot , y_{N,k_{r_{n}}}\in \rho_{N}^{k_{r_{n}}}x \cap V_{a_{n}}.$ If $k\notin \{k_{1},\cdot \cdot \cdot, k_{r_{n}}\},$ then we take elements $y_{1,k}\in \rho_{1}^{k}x,\cdot \cdot \cdot,  y_{N,k}\in \rho_{N}^{k}x$ arbitrarily (we know that such elements exist because $x\in \bigcap_{1\leq j \leq N} D_{\infty}(\rho_{j})$). With the sequence  
$y_{j,k}\in \rho_{j,k}(x)$ ($1\leq j\leq N$), the requirements of Definition \ref{4-skins-MLO-okay-novio}(i) satisfied.
\endproof

We are returning to the case of general case 
${\mathcal F}\in P(P({\mathbb N}))$ and ${\mathcal F}\neq \emptyset .$
If $\emptyset \notin {\mathcal F},$ then the ${\mathcal F}$-hypercyclicity (${\mathcal F}$-topological transitivity) of $G$ implies that $G$ is connected. On the other hand, if $\emptyset \in {\mathcal F}$ and $G$ is not connected, then
$G$ cannot be ${\mathcal F}$-hypercyclic (${\mathcal F}$-topologically transitive). In the sequel, we will consider only the case that $E(G)\neq \emptyset,$ when we clearly have ${\mathbb N} \in {\mathcal F};$ this will be our standing assumption in the sequel of this subsection. Since the associated binary relation $\rho$ is symmetric, we will also assume that
for each $A\in {\mathcal F}$ the assumption $A\neq \emptyset$ implies $A+2{\mathbb N}\subseteq A.$
Observe also that we have $S(U,V)=S(V,U)$ for any open non-empty subsets $U$ and $V$ of $G.$

As the next illustrative example shows, the notion of ${\mathcal F}$-hypercyclicity for simple graphs is far from being clear and easy understandable (see also \cite{extensions}):

\begin{example}\label{grafovi}
\begin{itemize}
\item[(i)]
Let $G=\{x_{1},x_{2},x_{3},x_{4}\}$ be equipped with discrete topology, let $G$ be the unoriented square $x_{1}x_{2}x_{3}x_{4},$ and let ${\mathcal F}$ be the collection of all non-empty subsets of ${\mathbb N}$ containing only odd elements. Then, for every $i\in {\mathbb N}_{4}$ and $n\in 2{\mathbb N}+1,$ we have that $x_{i}\notin \rho^{n}x_{i},$ which simply implies that the corresponding symmetric relation $\rho$ cannot possess any of the  introduced ${\mathcal F}$-dynamical properties from Definition \ref{4-skins-MLO-okay}; furthermore, $G$ is ${\mathcal F}$-hypercyclic (${\mathcal F}$-topologically transitive) iff $\{ 2{\mathbb N}, 2{\mathbb N}+1\} \subseteq {\mathcal F}$.    
\item[(ii)] Let the complete graph $K_{n}$ be equipped with discrete topology. Then the following holds:
\begin{itemize}
\item[(a)] $n=2:$ $K_{n}$ is ${\mathcal F}$-hypercyclic (${\mathcal F}$-topologically transitive) iff $\{ 2{\mathbb N}, 2{\mathbb N}+1\} \subseteq {\mathcal F}$.
\item[(b)]$ n\geq 3:$ $K_{n}$ is ${\mathcal F}$-hypercyclic (${\mathcal F}$-topologically transitive) iff ${\mathbb N} \setminus \{1\} \in {\mathcal F}.$
\end{itemize}
\end{itemize}
\item[(iii)] Let the complete graph $K_{n}$ be equipped with discrete topology. Then the following holds:
\begin{itemize}
\item[(a)] $n=2:$ The graphs
$K_{n},\cdot \cdot \cdot , K_{n},$ totally counted $N$ times, are ${\mathcal F}$-hypercyclic (${\mathcal F}$-topologically transitive) iff $\{ \emptyset , 2{\mathbb N}, 2{\mathbb N}+1\} \subseteq {\mathcal F}.$
\item[(b)] $n\geq 3:$ The graphs 
$K_{n},\cdot \cdot \cdot , K_{n},$ totally counted $N$ times, are ${\mathcal F}$-hypercyclic (${\mathcal F}$-topologically transitive) iff ${\mathbb N} \setminus \{1\} \in {\mathcal F}.$ 
\end{itemize}
\end{example}

Assume that $G$ is equipped with topology $\tau$ on $X$. Set
$$
{\mathrm S}_{G,\tau}:=\text{card}\bigl(\{S(U,V) : \emptyset \neq U, \ \emptyset \neq V,\ U,\ V\subseteq G\}\bigr).
$$
If $\tau$ is discrete topology on $X,$ then we simply write $
{\mathrm S}_{G}$ in place of $
{\mathrm S}_{G,\tau}.$

For connected bipartite graphs, the following result holds true:

\begin{thm}\label{moguce}
Let $G$ be a connected bipartite graph. Then, for every two open non-empty subsets $U$ and $V$ of $X,$ we have
\begin{align}\label{raw}
S(U,V)=L(U,V):=\bigl\{  d(u,v)+2k  : u\in U, v\in V,k\in {\mathbb N}_{0}\bigr\} \cap {\mathbb N}
\end{align}
and
\begin{align}\label{rah}
{\mathrm S}_{G,\tau} \leq {\mathrm S}_{G}\leq d(G) +\frac{(d(G))^{2}}{4}-\frac{1}{4}\chi_{2{\mathbb N}+1}(d(G)).
\end{align}
\end{thm}

\proof
Let $U$ and $V$ be given.  If $n\in S(U,V),$ then there exist nodes $u\in U,$ $v\in V$ and a walk in $G$ of length $n\in {\mathbb N}$ connecting $u$ and $v.$ If $d(u,v)=n,$ then clearly $n\in L(U,V);$ otherwise, $d(u,v)<n$ and $n-d(u,v)\in 2{\mathbb N},$ due to the fact that $G$ is bipartite, and we again have $n\in L(U,V).$ Conversely, if $n\in L(U,V),$ then there exist two nodes $u\in U,$ $v\in V$ and a number $k\in {\mathbb N}_{0}$ such that $n=d(u,v)+2k.$ If $u=v,$ then $k\in {\mathbb N}$ and there exists a walk in $G$ of length $2k$ connecting $u$ and $v=u$ because $G$ is connected. Otherwise, $u\neq v$  and there exists a walk in $G$ of length $d(u,v)$ connecting $u$ and $v.$ By the connectivity of $G,$ there exists a walk in $G$ of length $d(u,v)+2k=n,$ so that $n\in S(U,V)$
and \eqref{raw} holds. This implies that $S(U,V)$ is equal to some of the sets $1+2{\mathbb N}_{0},\cdot \cdot \cdot, d(G)+2{\mathbb N}_{0}$ or some of their finite non-empty unions. Since the inequality ${\mathrm S}_{G,\tau} \leq {\mathrm S}_{G}$ is trivial, for the proof of \eqref{rah}, it suffices to prove that there exist at most $\frac{(d(G))^{2}}{4}-\frac{1}{4}\chi_{2{\mathbb N}+1}(d(G))$ different finite non-empty unions of the sets $1+2{\mathbb N}_{0},\cdot \cdot \cdot, d(G)+2{\mathbb N}_{0}$. But, this simply follows from the fact that any such a union is of the form $A_{i,j}:=\{\min(i,j),\min(i,j)+2,\cdot \cdot \cdot, \max(i,j)-1\} \cup \{s\in {\mathbb N} : s\geq \max(i,j)\},$ where $1\leq i,j\leq d(G),$ $i$ is even, $j$ is odd, the fact that any such two sets $A_{i,j}$ and $A_{i',j'}$ differs if $(i, j)\neq (i',\ j'),$ $1\leq i,i',j,j'\leq d(G),$ $i,\ i'$ are even, $j,\ j'$ are odd, and the product principle.  
The proof of the theorem is thereby complete. 
\endproof

For the proof of inclusion $L(U,V) \subseteq S(U,V),$ we have not used the assumption that $G$ is bipartite. Therefore, we have:

\begin{prop}\label{mogucedq}
Let $G$ be a connected graph. Then, for every two  non-empty subsets $U$ and $V$ of $X,$ we have
$L(U,V) \subseteq S(U,V).$
\end{prop}

As Example \ref{grafovi}(i) shows, the estimate \eqref{raw} cannot be improved for connected bipartite graphs having four nodes (let us recall that ${\mathbb N}\in {\mathcal F}$ is our standing assumption). The situation is quite similar in general case because for the path $P_{n},$ where $n\geq 2,$ we have 
$$
d(P_{n})=n-1 \mbox{ and } 
S_{G}=n-1+\frac{(n-1)^{2}}{4}-\frac{1}{4}\chi_{2{\mathbb N}+1}(n-1).
$$

By $\vartheta(G)$ we denote the smallest number, if such exists, satisfying that any two nodes $u$ and $v$ of $G$ can be connected by an even walk of length $\leq \vartheta(G)$ and an odd walk of length $\leq \vartheta(G);$ otherwise, we set $\vartheta(G):=+\infty.$ 

\begin{thm}\label{reza}
Let $G$ be a connected graph.
Then we have
\begin{align}\label{rahmanji}
{\mathrm S}_{G} \leq \Biggl \lfloor \frac{1}{4}\Bigl( \vartheta^{2}(G)+2\vartheta(G)+1  \Bigr) \Biggr \rfloor.
\end{align}
\end{thm}

\proof
Clearly, it suffices to examine the case in which $\vartheta(G)<+\infty.$ In this case, we have
\begin{align}\label{rezonovanje}
B:=\{ s\in {\mathbb N} : s\geq \vartheta(G) \} \subseteq S(U,V)
\end{align}
for any two non-empty subsets $U$ and $V$ of $G.$ Set $m_{U,V}:=\min (S(U,V)).$ If $m_{U,V}=\vartheta(G),$ then $S(U,V)=\{ s\in {\mathbb N} : s\geq \vartheta(G) \}.$  If $1\leq i=m_{U,V}<\vartheta(G),$ then $S(U,V) =(i+2{\mathbb N}_{0}) \cup B$ or there exists a natural number  
$i'\in (i, \vartheta(G))$ such that
$i'-i$ is an odd number and
$S(U,V)=\{s\in {\mathbb N} : s\in i+2{\mathbb N}_{0},\ s<i'\} \cup \{s \in {\mathbb N} : s\geq i'\}.$ Hence, if
$\vartheta(G)-i$ is an even number, we have at most $(\vartheta(G)-i+2)/2$ different possibilities for $S(U,V),$ while if 
$\vartheta(G)-i$ is an odd number, we have at most $(\vartheta(G)-i+1)/2$ different possibilities for $S(U,V).$ Summa summarum, 
\begin{align}\label{nekt}
S_{G} \leq 1+\sum_{i=1}^{\vartheta(G)-1}\frac{\vartheta(G)-i+2}{2}-\sum_{i\in I}\frac{1}{2},
\end{align}
where $I=\{ i \in {\mathbb N}_{\vartheta(G)-1} : \vartheta(G)-1-i \mbox{ is odd}\}.$
Then the estimate \eqref{rahmanji} follows from a simple computation involving \eqref{nekt} and the equality $\lfloor \frac{1}{4}( \vartheta^{2}(G)+2\vartheta(G)+1 )  \rfloor=\lfloor \frac{1}{4}( \vartheta^{2}(G)+2\vartheta(G) )  \rfloor ,$ holding for even numbers $\vartheta(G)$.
\endproof

If a connected graph $G$ is not bipartite, then it contains a closed contour of odd length as a subgraph. If $l\in {\mathbb N}$ is any number such that the closed contour $C_{2l+1}$ is a subgraph of $G,$ then it is very elementary to prove that $\vartheta (G) \leq \max_{u,v\in G}[2d(u,C_{2l+1})+d(u,v)+2l+1],$ where $d(u,C_{2l+1}):=\inf\{ d(u,w) : w\in C_{2l+1} \}.$ 

\begin{rem}\label{dragosi}
For a non-bipartite connected graph $G$, the estimate \eqref{rezonovanje} follows immediately from the facts that the index of $G$ is strictly greater than $1,$ any component of principal eigenvector of $G$ is strictly positive and an application of \cite[Theorem 2.2.5]{dragos}; see \cite{dragos} for the notion.
\end{rem}

Let $m,\ n\in {\mathbb N},$ let $X=\{x_{1},x_{2},\cdot \cdot \cdot , x_{n}\},$ $Y=\{y_{1},y_{2},\cdot \cdot \cdot , y_{m}\},$ and let $G_{i}(X,Y)$ be a bipartite graph with colored classes $X$ and $Y$ ($1\leq i \leq N$). Suppose that $X\cup Y$ is equipped with discrete topology and ${\mathcal F}=P(P({\mathbb N})) \setminus \{\emptyset \}.$
As indicated in \cite{marina-prim}, the graphs $G_{1}(X,Y),G_{2}(X,Y),\cdot \cdot \cdot,G_{N}(X,Y)$ cannot be strongly $d{\mathcal F}$-hypercyclic (strongly $d{\mathcal F}$--topologically transitive); by Proposition \ref{poka}, it readily follows that $G_{1}(X,Y),G_{2}(X,Y),\cdot \cdot \cdot,G_{N}(X,Y)$ cannot be $d{\mathcal F}$-hypercyclic ($d{\mathcal F}$--topologically transitive). On the other hand,
by our comment from Remark \ref{dragosi},
for arbitrary non-bipartite connected graphs $G_{1}, G_{2},\cdot \cdot \cdot, G_{N}$ and for arbitrary non-empty subsets $V_{1},V_{2},\cdot \cdot \cdot,V_{N}$ of $ {\mathbb N}_{n},$ we always have the existence of a positive integer $k_{0}\in {\mathbb N}$ such that
$$
\bigl[k_{0},\infty\bigr) \cap {\mathbb N} \subseteq \bigl\{ k\in {\mathbb N} : (\forall i\in U)(\forall s\in {\mathbb N}_{N})(\forall j_{s}\in V_{s}) a_{ij_{s}}^{k;s}\geq 1\bigr\}.
$$
Taking into account Proposition \ref{boze-prim}, the above immediately implies:

\begin{thm}\label{radio}
Let $G_{1}, G_{2},\cdot \cdot \cdot, G_{N}$ be non-bipartite connected graphs, and let ${\mathcal F}=P(P({\mathbb N})) \setminus \{\emptyset \}.$ Then $G_{1}, G_{2},\cdot \cdot \cdot, G_{N}$ are always (strongly) $d{\mathcal F}$-hypercyclic (strongly $d{\mathcal F}$-topologically transitive) and any element of $X$ is a (strong) $d{\mathcal F}$-hypercyclic vector of $G_{1}, G_{2},\cdot \cdot \cdot, G_{N}.$
\end{thm}

With the notion introduced in \cite{marina-prim}, we also have that $G_{1}, G_{2},\cdot \cdot \cdot, G_{N}$ are $d$-topologically mixing.

We would like to propose the following problem:\vspace{0.1cm}

\noindent {\sc Problem 1.} Given a connected non-bipartite graph $G$ with $n\geq 2$ nodes, find as better as possible upper bound for ${\mathrm S}_{G}$ in terms of $d(G)$ and $n.$ \vspace{0.1cm}

Theorem \ref{moguce} and Proposition \ref{reza} can be reconsidered for disjointness; Problem 1 can be reformulated in this context, as well. For the sake of brevity, we will skip all related details about these questions. 

\subsection{Results for digraphs and tournaments}\label{digraphs} 

A digraph is any pair $(G,\rho),$ where $G$ is a finite non-empty set
and $\rho \subseteq (G\times G) \setminus \Delta_{G};$ hence, in our definition, we do not allow $G$ to contain any loop. We will consider only finite non-trivial digraphs henceforth ($|G|>1$). The
elements in $G$ and $\rho$ are called points (vertices, nodes) and arcs respectively; if arc
$(x, y)\in \rho,$ then we say that $x$ is adjacent to $y$ and write 
$xy$ for arc$ (x, y).$ Two vertices $x$ and $y$ of a digraph $G$
are said to be nonadjacent iff $(x, y)\notin \rho$ and $(y,x)\notin \rho.$ If we replace each  arc$ (x, y)$ in $G$ by symmetric pairs $(x,y)$ and $(y,x)$ of arcs, we obtain the underlying simple graph ${\mathrm G}$ associated to $G.$
The notions of outdegree $d^{+}(x),$ indegree $d^{-}(x)$ and degree $d(x):=d^{+}(x)+d^{-}(x)$ of a vertex $x\in G$ as well as the notions of Hamiltonicity of $G,$ a semi-walk in $G,$ a walk in $G$ and their lengths are defined usually (\cite{booksa}). A digraph
$G$ is called asymmetric iff $\rho$ is an anti-symmetric relation. If $G$ is asymmetric digraph and $G$ is ${\mathcal F}$-hypercyclic (${\mathcal F}$-topologically transitive), then 
for each set $A\in {\mathcal F}$ we have $2\notin A.$
Let us recall that a tournament $T$ is a digraph in which any two different nodes are connected by exactly one arc. The set of nodes of any digraph $G$ (tournament $T$) considered below will be $X=V(G)=\{x_{1},x_{2},\cdot \cdot \cdot, x_{n}\}.$

Let us recall that a digraph $(G,\rho)$ is said to be strongly connected iff for any two different points $x$ and $y$ from $G$ there is an oriented $x-y$ walk, while
$(G,\rho)$ is said to be weakly connected iff for any two different points $x$ and $y$ from $G$ there is an $x-y$ semi-walk, which is equivalent to say that the underlying simple graph ${\mathrm G}$ associated to  $G$
is connected (\cite{booksa}). For various generalizations, see \cite{marina-prim}.
By $[A(G)]_{1\leq i,j\leq n}$ we denote the adjacency matrix of $G.$ 

For any digraph $G$ (digraphs $G_{1},\cdot \cdot \cdot, G_{N}$), denote by ${\mathrm G}$ (${\mathrm G}_{1},\cdot \cdot \cdot, {\mathrm G}_{N}$) the associated simple graphs defined as above. The notions introduced in Definition \ref{4-skins-MLO-okay} (Definition \ref{4-skins-MLO-okay-novio} and Definition \ref{4-skins-MLO-okay-novio'frim}) can be used to define 
${\mathcal F}_{w}$-hypercyclic and and ${\mathcal F}_{w}$-topologically transitive properties ($d{\mathcal F}_{w}$-hypercyclic and $d{\mathcal F}_{w}$-topologically transitive properties) of $G$ ($G_{1},\cdot \cdot \cdot, G_{N}$). For example, we say that $G$ is strongly ${\mathcal F}_{w}$-hypercyclic
iff the associated simple graph ${\mathrm G}$ is strongly ${\mathcal F}$-hypercyclic, while $G_{1},\cdot \cdot \cdot, G_{N}$ are said to be $d{\mathcal F}_{w}$-topologically transitive
iff ${\mathrm G}_{1},\cdot \cdot \cdot, {\mathrm G}_{N}$ are $d{\mathcal F}$-topologically transitive, and so on and so forth. In such a way, we extend the notion of 
$d$-weakly connected digraphs (\cite{marina-prim}). It is clear that any ${\mathcal F}$-hypercyclic property ($d{\mathcal F}$-hypercyclic property) implies the corresponding ${\mathcal F}_{w}$-hypercyclic property ($d{\mathcal F}_{w}$-hypercyclic property). The same holds for topological transitivity.

The numbers ${\mathrm S}_{G,\tau}$ and ${\mathrm S}_{G}$ are meaningful for digraphs, as well, but calculating upper bounds for the number ${\mathrm S}_{G}$ is not so simple task for digraphs. Let us only note that in a primitive digraph $G$ (this means that there is a positive integer $k\in {\mathbb N}$ such that there is a walk of length
$k$ from each vertex $u$ to each vertex $v$  (possibly $u$ again) of $G;$ the smallest
integer $k$ with this property is said to be the exponent of $G$), for any two open non-empty subsets $U$ and $V$ of $X,$ we have $\{k\in {\mathbb N} :k\geq exp(G)\}\subseteq S(U,V),$
where $exp(G)$ denotes the exponent of $G.$ Unless stated otherwise, we assume henceforth that ${\mathcal F}=P(P({\mathbb N})) \setminus \{\emptyset \}.$ 

The following result is closely connected with Proposition \ref{poka}:

\begin{prop}\label{vaterpolo}
\begin{itemize}
\item[(i)]
If $G$ is equipped with discrete topology, then $G$ is ${\mathcal F}$-hypercyclic (${\mathcal F}$-topologically transitive) iff $G$ is strongly ${\mathcal F}$-hypercyclic (strongly ${\mathcal F}$-topologically transitive).
\item[(ii)] Suppose that $G_{1},G_{2},\cdot \cdot \cdot, G_{N}$ are given digraphs and $G_{i}$ is equipped with discrete topology on $X$ ($1\leq i\leq N$). Then $G_{1},G_{2},\cdot \cdot \cdot, G_{N}$ are $d{\mathcal F}$-hypercyclic ($d{\mathcal F}$-topologically transitive) iff $G_{1},G_{2},\cdot \cdot \cdot, G_{N}$ are strongly $d{\mathcal F}$-hypercyclic (strongly $d{\mathcal F}$-topologically transitive).
\end{itemize}
\end{prop}

\proof
We present only the main points of the proof of (ii) for $d{\mathcal F}$-hypercyclicity and strong $d{\mathcal F}$-hypercyclicity, which is very similar to that of Proposition \ref{poka}. It suffices to show that any $d{\mathcal F}$-hypercyclic vector of $G_{1},G_{2},\cdot \cdot \cdot, G_{N}$ is likewise a strong $d{\mathcal F}$-hypercyclic vector of $G_{1},G_{2},\cdot \cdot \cdot, G_{N}.$ To see this, we can copy the arguments given for simple graphs because we can always find a strictly increasing sequence $k_{1}<k_{2}<\cdot \cdot \cdot<k_{r_{n}}$ satisfying the properties stated in the proof of Proposition \ref{poka}, due to our assumption that any digraph $G_{i}$ is equipped with discrete topology on $X$ ($1\leq i\leq N$) and the fact that for each $d{\mathcal F}$-hypercyclic vector $x$ of $G_{1},G_{2},\cdot \cdot \cdot, G_{N}$ there exists a positive integer $k\in {\mathbb N}$ such that, for every $j\in {\mathbb N}_{N},$ there exists an $(x-x)$ walk in $G_{j}$ of length $k.$ 
\endproof

Using Proposition \ref{vaterpolo}, we can rephrase a great number of our results established in \cite{marina-prim} for strong 
${\mathcal F}$-hypercyclicity (strong ${\mathcal F}$-topological transitivity) and strong $d{\mathcal F}$-hypercyclicity (strong $d{\mathcal F}$-topological transitivity). For example, we have the following:
\begin{itemize}
\item[(i)]
Let $G$ be a tournament equipped with discrete topology. Then $G$ is strongly ${\mathcal F}$-hypercyclic iff the indegree of any vertex is strictly positive.
\item[(ii)] Let $G$ be a digraph equipped with discrete topology, satisfying that for any two vertices $x,\ y$ in $G$ such that 
$xy$ is not an arc in $G$ one has $d^{+}(x)+d^{-}(y)\geq n-1.$  Then $\rho$ is strongly ${\mathcal F}$-hypercyclic iff the indegree of any vertex is strictly positive.
\item[(iii)] Let $n\geq 4,$ and let $T_{1},\ T_{2},\cdot \cdot \cdot,\ T_{N}$ be tournaments equipped with discrete topologies. Then $T_{1},\ T_{2},\cdot \cdot \cdot,\ T_{N}$ are strongly $d{\mathcal F}$-topologically transitive iff
$T_{j}$ is strongly connected for all $j\in {\mathbb N}_{N}$ iff
$T_{j}$ is a Hamiltonian tournament for all $j\in {\mathbb N}_{N}$.
\end{itemize}

It is worth noting that Proposition \ref{vaterpolo} does not hold if $G$ (some of $G_{i}'$s for $1\leq i\leq N$) is equipped with topology that is not discrete:

\begin{example}\label{pende}
\begin{itemize}
\item[(i)]
Let $n=5$ and $G$ be a digraph with the associated binary relation $\rho=\{(x_{3},x_{2}), (x_{3},x_{5}), (x_{2},x_{1}), (x_{1},x_{4}), (x_{4},x_{1})\}.$ Suppose that $G$ is equipped with topology $\tau=\{\emptyset, \{x_{2}\}, \{x_{5}\},\{x_{2},x_{5}\}, \{x_{1},x_{2},x_{3},x_{4},x_{5}\}\}.$ Then $x_{3}$ is the only ${\mathcal F}$-hypercyclic vector of $G$ and therefore $G$ is ${\mathcal F}$-hypercyclic. On the other hand, $G$ is not ${\mathcal F}$-topologically transitive because there is no $(x_{2}-x_{5})$ walk in $G.$ In this concrete example, $G$ is not strongly ${\mathcal F}$-hypercyclic because $x_{3}$ is not a strong ${\mathcal F}$-hypercyclic vector of $G;$ this follows from the fact that for the sequence $(y_{n})$ satisfying the requirements from Definition \ref{4-skins-MLO-okay}(i) we need to have $y_{1}=x_{2}$ and $y_{1}=x_{5},$ which is a contradiction. If we consider $N$ copies of digraph $G,$ then it can be simply verified that the obtained tuple is $d{\mathcal F}$-hypercyclic, not strongly $d{\mathcal F}$-hypercyclic and not $d{\mathcal F}$-topologically transitive.
\item[(ii)] Let $n=4$ and $G$ be a digraph with the associated binary relation $\rho=\{(x_{1},x_{2}), (x_{1},x_{3}), (x_{3},x_{4}), (x_{4},x_{3})\}.$ Suppose that $G$ is equipped with topology $\tau=\{\emptyset, \{x_{2}\}, \{x_{3}\},\{x_{2},x_{3}\}, \{x_{1},x_{2},x_{3},x_{4}\}\}.$ 
Consider the $N$ copies of digraph $G;$
then $x_{1}$ is the only $d{\mathcal F}$-hypercyclic vector of obtained tuple ${\mathrm T},$ which is $d{\mathcal F}$-hypercyclic and not strongly $d{\mathcal F}$-hypercyclic. Otherwise, $x_{1}$ needs to be a strong $d{\mathcal F}$-hypercyclic vector of ${\mathrm T},$ which is a contradiction because the choices $V_{1}=V_{2}=\{x_{2}\}$ and $V_{1}=\{x_{2}\},\ V_{2}=\{x_{3}\}$ impose that, for the sequence $(y_{j,k})_{1\leq j\leq N,k\in {\mathbb N}}$ in Definition \ref{4-skins-MLO-okay-novio}(i), we must have $y_{2,1}=x_{2}$ and $y_{2,1}=x_{3}.$ Observe also that $G,$ resp. ${\mathrm T},$ is not ${\mathcal F}$-topologically transitive, resp. $d{\mathcal F}$-topologically transitive. 
\end{itemize}
\end{example}

The values $n= 5$ and $n=4$ are optimal, as the next two propositions indicate:

\begin{prop}\label{pende-primp}
Let $n\leq 4$ and let $G$ be a given digraph. Then $G$ is ${\mathcal F}$-hypercyclic iff $G$ is strongly ${\mathcal F}$-hypercyclic. 
\end{prop}

\proof
We will consider only the most complicated case $n=4;$ the proof for ${\mathcal F}$-hypercyclicity goes as follows. Let $x=x_{1}$ be a ${\mathcal F}$-hypercyclic vector for $G$ and let  $V_{1},\cdot \cdot \cdot, V_{m}$ be all  non-empty open subsets of $X.$ Then there exist a natural number $k_{i}$ and an element $y_{k_{i}}\in \rho^{k_{i}}x \cap V_{i}$ ($1\leq i\leq m$). If there exists an $(x_{1}-x_{1})$ walk in $X,$ it is clear that $k_{1},\cdot \cdot \cdot, k_{i}$ can be chosen arbitrarily large and it is trivial to show that, in this case, $x_{1}$ needs to be a strong ${\mathcal F}$-hypercyclic vector for $G.$ Suppose that there is no $(x_{1}-x_{1})$ walk in $X$ and $x_{1}\neq y_{k_{i}}$ for $1\leq i\leq m;$ since $x\in D_{\infty}(\rho),$ it readily follows that the unoriented segment $x_{2}x_{3},$ $x_{2}x_{4}$ or $x_{3}x_{4}$ belong to $\rho,$ with the meaning clear, or a closed contour connecting $x_{2},$ $x_{3}$ and $x_{4}$ belongs to $\rho,$ again with the meaning clear. The last case is trivial because we can reach the points $x_{2},$ $x_{3}$ and $x_{4}$ by  walks of arbitrarily large length, starting from $x_{1},$ so that $x_{1}$ needs to be a strong ${\mathcal F}$-hypercyclic vector for $G$ by an elementary line of reasoning. Otherwise, we may assume that the unoriented segment $x_{3}x_{4}$ belongs to $\rho.$ If $x_{3}x_{2}$ or $x_{4}x_{2}$ is an arc in $G,$ then it is clear that we can reach the points $x_{2},$ $x_{3}$ and $x_{4}$ by walks of arbitrarily large lengths, starting from $x_{1},$ so that the proof is complete. If this is not the case, then $x_{1}x_{2}$ may or may not be an arc in $G.$ In the first case, if some of elements $y_{1},\cdot \cdot \cdot,y_{m}$ is equal to $x_{2},$ then we may assume without of generality that at most one of these elements equal to $x_{2}$  
(because the set $\{x_{2}\}$ can be listed in the sequence $V_{1},\cdot \cdot \cdot,V_{m}$ at most once and 
the
points $x_{3}$ and $x_{4}$ can be reached by walks from $x_{1}$). In the second case, any of elements $y_{1},\cdot \cdot \cdot,y_{m}$ cannot be equal to $x_{2}$ and any of numbers $k_{1},\cdot \cdot \cdot,k_{m}$ can be chosen arbitrarily large because the points $x_{3}$ and $x_{4}$ can be reached by walks of arbitrarily large lengths, starting from $x_{1}.$ The proof of proposition is completed.
\endproof

\begin{prop}\label{pende-primp=disjoint}
Let $n\leq 3$ and let $G_{1},G_{2},\cdot \cdot \cdot, G_{N}$ be given digraphs. Then $G_{1},G_{2},\cdot \cdot \cdot, G_{N}$ are $d{\mathcal F}$-hypercyclic iff $G_{1},G_{2},\cdot \cdot \cdot, G_{N}$ are strongly $d{\mathcal F}$-hypercyclic. 
\end{prop}

\proof
The proof is very similar to that of Proposition \ref{poka}.
Let $x_{1}$ denote a $d{\mathcal F}$-hypercyclic vector of $G_{1},G_{2},\cdot \cdot \cdot, G_{N}.$ Then some of the unoriented arcs $x_{1}x_{2},$ $x_{1}x_{3},$ $x_{2}x_{3}$ or a closed oriented contour of length $3$  is contained in any $G_{i}$ ($1\leq i\leq N$); this implies that, in the proof of Proposition \ref{poka}, we can choose $k_{2}>k_{1},$ because the number $k_{2}$ can be replaced therein, optionally, with any number $k_{2}+6s,$ where $s\in {\mathbb N}.$ Keeping this in mind, we can repeat literally the arguments given in the proof of Proposition \ref{poka}.
\endproof

\begin{rem}\label{prcojed}
In our previous analyses, we have constructed a tournament (two tournaments) having two nodes and equipped with certain topology (topologies) that 
is ${\mathcal F}$-topologically transitive but not strongly ${\mathcal F}$-topologically transitive ($d{\mathcal F}$-topologically transitive but not strongly $d{\mathcal F}$-topologically transitive). Hence, the statements of Proposition \ref{pende-primp} and Proposition \ref{pende-primp=disjoint} do not hold for ${\mathcal F}$-topological transitivity and strong ${\mathcal F}$-topological transitivity.
\end{rem}

Any tournament ($N$ tournaments) having two nodes cannot be ${\mathcal F}$-hypercyclic because the domain of square of the associated binary relation (relations) is the empty set. In the case that $n\geq 3,$ we have the following result closely connected with Proposition \ref{pende-primp} and Proposition \ref{pende-primp=disjoint}:

\begin{thm}\label{pende-bn}
\begin{itemize}
\item[(i)] Suppose that $n\geq 3$ and $T$ is a given tournament. Then $T$ is ${\mathcal F}$-hypercyclic iff $T$ is strongly ${\mathcal F}$-hypercyclic.
\item[(ii)] Suppose that $n\leq 4$ and $T_{1},T_{2},\cdot \cdot \cdot, T_{N}$ are given tournaments. Then $T_{1},T_{2},\cdot \cdot \cdot, T_{N}$ are $d{\mathcal F}$-hypercyclic iff $T_{1},T_{2},\cdot \cdot \cdot, T_{N}$ are strongly $d{\mathcal F}$-hypercyclic. 
\item[(iii)] Suppose that $n\geq 5.$ Then we can always find tournaments $T_{1},T_{2},\cdot \cdot \cdot, T_{N}$ that are $d{\mathcal F}$-hypercyclic and not strongly $d{\mathcal F}$-hypercyclic. 
\end{itemize}
\end{thm}

\proof
We will prove (i) by induction. If $n=3,$ then $T$ is isomorphic to a Hamiltonian contour, which is clearly strongly ${\mathcal F}$-hypercyclic, or to the tournament with the set of nodes $x_{1}x_{2},$ $x_{1}x_{3}$ and $x_{2}x_{3},$ which cannot be  ${\mathcal F}$-hypercyclic because the cube of associated binary relation has empty domain. Suppose that the statement of proposition holds for each tournament having strictly less than $n>3$ nodes and let us prove the statement for an arbitrary tournament $T$ having $n$ nodes. 
Let $x_{i_{1}}$ be an ${\mathcal F}$-hypercyclic vector of $T.$ Due to the famous theorem of L. R\' edei (see e.g. \cite{petrovic}), there exists a Hamiltonian path, say $x_{1}\mapsto x_{2} \mapsto x_{3} \mapsto \cdot \cdot \cdot 
\mapsto x_{n},$ in $T.$ If there exists adjacent nodes $x_{j}$ and $x_{l}$ for some $j\leq i_{1}$ and $l\geq i_{1},$ $j,\ l\in {\mathbb N}_{n},$ then there is a closed $(x_{i_{1}}-x_{i_{1}})$-walk in $T$ and the statement trivially holds. Otherwise, the set of nodes which can be reached from $x_{i_{1}}$ by a closed walk is $X'\equiv \{x_{i_{1}+1},\cdot \cdot \cdot, x_{n}\}.$ If $i_{1}>1,$ then the subtournament $T'$ induced by the set of nodes $X'$ is ${\mathcal F}$-hypercyclic with respect to the subspace topology on $X'.$ By induction hypothesis, $T'$ is strongly ${\mathcal F}$-hypercyclic which clearly implies that $T$ is strongly ${\mathcal F}$-hypercyclic, as well. It remains to be considered the case that $i_{1}=1.$ The result trivially follows if the initial topology is anti-discrete or the set $\{x_{1}\}$ is an open set of the initial topology, because then there exists a closed $(x_{1}-x_{1})$-walk in $T.$ If this is not the case,
set $y_{i}:=x_{i+1}$ for $1\leq i\leq n-1$ and take $y_{i}\in \rho^{i}x_{1}$ arbitrarily for $i\geq n.$ Then for any open set $V$ of the initial topology there exists an integer $i\in {\mathbb N}_{n} \setminus \{1\}$ such that $x_{i}\in V.$ It is clear that the set $\{k\in {\mathbb N} : y_{k}\in V\}$ is non-empty because it contains the element $x_{i-1}.$
This completes the proof of (i). The proof of (ii) for $n=2$ and $n=3$ is simple. If $n=4,$ 
assume that $x_{1}$ is a $d{\mathcal F}$-hypercyclic vector of $T_{1},T_{2},\cdot \cdot \cdot, T_{N}.$ Then Proposition \ref{komponente} and the fact that $x_{1}\in \bigcap_{1\leq j\leq N, n\in {\mathbb N}}D_{\infty}(\rho_{j})$ implies that a closed oriented contour of length $3$ or a closed oriented contour of length $4$  is contained in any $G_{i}$ ($1\leq i\leq N$). Then 
we can argue as in the proofs of Proposition \ref{poka} and Proposition \ref{pende-primp=disjoint} because we can construct the sequence \
$(k_{s})_{1\leq s\leq r_{n}}$ such that
$k_{1}<k_{2}<\cdot \cdot \cdot < k_{r_{n}-1}<k_{r_{n}};$ this follows from the fact that any point reachable from $x_{1}$ by a walk of length $k$ is also reachable from $x_{1}$ by a walk of length $k+12s,$ where $s\in {\mathbb N}.$ For a counterexample in (iii), let $T_{1}$ be any tournament containing the following sets of arcs $\{x_{1}x_{i} :  1\leq i\leq n\},$ $ \{x_{2}x_{i} : 3\leq i\leq n\},$ $ \{x_{3}x_{i} : i=4\mbox{ or }6\leq i\leq n\},$  $\{x_{4}x_{i} : 5\leq i \leq n\}$ and $\{x_{5}x_{i} : i=3\mbox{ or } 6\leq i\leq n\}.$ 
Let 
$T_{2}$ be any tournament containing the following sets of arcs $\{x_{1}x_{i} :  1\leq i\leq n\},$ $ \{x_{3}x_{i} : i=2\mbox{ or }4\leq i\leq n\},$ $ \{x_{2}x_{i} : 4\leq i\leq n\},$  $\{x_{4}x_{i} : 5\leq i \leq n\}$ and $\{x_{5}x_{i} : i=2\mbox{ or } 6\leq i\leq n\}.$ Let $\tau=\{\emptyset, \{x_{2}\},\{x_{3}\},\{x_{2},x_{3}\}, X\}.$ Then it can be easily seen that 
$x_{1}$ is a unique $d{\mathcal F}$-hypercyclic vector of $T_{1},T_{2},\cdot \cdot \cdot, T_{N}$ and that $T_{1},T_{2},\cdot \cdot \cdot, T_{N}$ are $d{\mathcal F}$-hypercyclic but not strongly $d{\mathcal F}$-hypercyclic because $x_{1}$ is not a strong $d{\mathcal F}$-hypercyclic vector of $T_{1},T_{2},\cdot \cdot \cdot, T_{N};$ this can be seen by plugging $V_{1}=V_{2}=\{x_{2}\}$ and $V_{1}=V_{2}=\{x_{3}\}$ in Definition \ref{4-skins-MLO-okay-novio}(i), which immediately forces that $y_{1,1}=x_{2}$ and $y_{1,1}=x_{3},$ a contradiction.
\endproof

The strong connectivity of a tournament $T$ is equivalent with its irreducibility (\cite{moon-prim}). Applying \cite[Theorem 1]{moon-prim}, we get that
the adjacency matrix $[A(T)]$ of $T$ is primitive, i.e., there exists a natural number $q\in {\mathbb N}$ such that any element of
the matrix $[A(T)]^{q}$ is strictly positive. Due to the strong connectivity of $T,$ we get that any element of the matrix $[A(T)]^{q'}$ is strictly positive for all $q'\geq q$ (see also \cite{marina-prim}). If we denote by $d$ and $e$ the diameter and exponent of $T$, then $d\leq e\leq d+3$ and we can take $e=q\leq d+3$ due to \cite[Theorem 19]{moon}; see \cite[Section 12, Section 13]{moon} for the notion and more details on the subject.
This immediately implies the following result:

\begin{prop}\label{idiot}
\begin{itemize}
\item[(i)] Let $T$ be a strongly connected tournament with exponent $e.$ If $T$ is ${\mathcal F}$-hypercyclic (${\mathcal F}$-topologically transitive), then for each set $A\in {\mathcal F}$ we have $\{ n\in {\mathbb N} : n\geq e \}\subseteq A.$ 
\item[(ii)] Let $T_{1},\cdot \cdot \cdot, T_{N}$ be strongly connected tournaments, and let $e$ be the largest value of their exponents.
If $T_{1},\cdot \cdot \cdot, T_{N}$ are $d{\mathcal F}$-hypercyclic ($d{\mathcal F}$-topologically transitive), then for each set $A\in {\mathcal F}$ we have $\{ n\in {\mathbb N} : n\geq e \}\subseteq A.$ 
\end{itemize}
\end{prop}

It is a well known fact that the exponent $e$ of a primitive tournament $T$ satisfies $3\leq e\leq n+2$ 
if $n\geq 5.$ Furthermore, if $n\geq 6,$
then
\cite[Theorem 20]{moon} implies the existence of a primitive tournament $T$ with exponent $e\in [3,n+2]$ given in advance. 
For such a tournament $T,$ we have the existence of two elements $x_{i}$ and $x_{j},$ for some $i,\ j \in {\mathbb N}_{n},$ such that there is no ($x_{i}-x_{j})$ walk in $T$ of length $e-1.$ This implies that $x_{i}$ cannot be an ${\mathcal F}$-hypercyclic vector of $T$ if each set $A\in {\mathcal F}$ contains the number $e-1.$ 

We close the paper by proposing some open problems and observations for tournaments having four vertices.
There exist four non-isomorphic tournaments with four vertices. The only one of them is transitive (acyclic) tournament and this tournament is ${\mathcal F}$-topologically transitive for discrete topology iff $\{\emptyset , \{1\}, \{2\}, \{1,2,3\}\} \subseteq {\mathcal F}.$ The second (third) one is obtained as the union of the closed contour $x_{2}x_{3}x_{4}$ and arcs $x_{1}x_{2},$ $x_{1}x_{3},$ $x_{1}x_{4}$ ($x_{2}x_{1},$ $x_{3}x_{1},$ $x_{4}x_{1}$) and these tournaments are ${\mathcal F}$-topologically transitive for discrete topology iff ${\mathcal F}$ contains the sets $\emptyset ,$ $ 3{\mathbb N},$ $ 1+3{\mathbb N},$
$2+3{\mathbb N}$ and their finite unions. The fourth one is obtained as the union of closed contour $x_{1}x_{2}x_{3}x_{4}$ and arcs $x_{1}x_{3},$ $x_{2}x_{4}.$ Direct computation of powers of corresponding adjacency matrix shows that this tournament is ${\mathcal F}$-topologically transitive for discrete topology iff ${\mathcal F}$ contains the sets ${\mathbb N} \setminus \{2\},$
${\mathbb N} \setminus \{1,4\},$ ${\mathbb N} \setminus \{2,3,6\},$ ${\mathbb N} \setminus \{1,2,3,6\},$ ${\mathbb N} \setminus \{1,3,4,7\},$ ${\mathbb N} \setminus \{1,2,4,5,8\}$ and their finite unions. The obtained results seem to be very dissociated and, because of that, we would like to propose the following problem: \vspace{0.1cm}

{\sc Problem 2.} Let  $n\geq 2,$ and let $A_{i}\subseteq {\mathbb N}$ be a given set ($1\leq i \leq (2^{n}-1)(2^{n}-1)$). Find necessary and sufficient conditions for the existence of a tournament $T_{n}$ with $n$ vertices such that $\{S(U,V) : \emptyset \neq U, \ \emptyset \neq V,\ U,\ V\subseteq V(T)\}=\{ A_{i} : 1\leq i\leq	 (2^{n}-1)(2^{n}-1) \}.$
\vspace{0.1cm}

It is also meaningful to ask the following:\vspace{0.1cm}

{\sc Problem 3.} Set $a_{n}:=\{ {\mathrm S}_{T_{n}} : T_{n} \mbox{ is a tournament with } n \mbox{ vertices}\},$ $n\geq 2.$ Find some upper bounds for $a_{n}$ ($n\geq 2$) and an asymptotic behaviour of the sequence $(a_{n})_{n\geq 2}.$ 
\vspace{0.1cm}

{\sc Problem 4.} Let  $n\geq 2.$ Construct a tournament $T_{n}$ with $n$ vertices such that ${\mathrm S}_{T_{n}}=a_{n}.$ How many non-isomorphic tournaments with $n$ vertices satisfy this equality? \vspace{0.1cm}

{\sc Problem 5.} Reconsider Problem 2, Problem 3 and Problem 4 for strongly connected  tournaments and some classes of asymmetric digraphs.

{\small

\end{document}